\documentclass[titlepagea4paper,11pt]{amsart}
\usepackage[arrow,matrix]{xy}
\usepackage{amssymb}
\usepackage{amsbsy,amsmath}
\newtheorem{thm}{Theorem}[section]
\newtheorem{prop}[thm]{Proposition}
\newtheorem{lemma}[thm]{Lemma}

\newcommand{\seq}{\longrightarrow}

\newcommand{\dd}{\partial\overline\partial}

\newcommand{\Ric}{\mbox{Ric}}

\newcommand{\bb}{\mathbb}
\newcommand{\4}[4]{_{#1\bar#2#3\bar#4}}
\newcommand{\ii}[2]{_{#1\bar#2}}
\newcommand{\II}[2]{^{#1\bar#2}}

\newtheorem{example}[thm]{Example}
\newtheorem{construction}[thm]{Construction}

\title{On asymptotics of complete Ricci-flat K\"ahler metrics on open manifolds}
\author{Bert Koehler and Marco K\"uhnel}
\subjclass[2000]{32W20; 32Q15}
\keywords{Ricci-flat, K\"ahler metrics, Monge-Amp\`ere equation, maximum principle}
\address{Bert Koehler\\Debeka Hauptverwaltung\\Abteilung LV/M\\Ferdinand-Sauerbruch-Str. 18\\56058 Koblenz, Germany}
\email{bert.koehler@debeka.de}
\address{Marco K\"uhnel\\Otto-von-Guericke-Universit\"at Magdeburg\\FMA / IAN\\Postfach 4120\\
39016 Magdeburg\\Germany}
\email{marco.kuehnel@mathematik.uni-magdeburg.de}
\thanks{The authors acknowledge gratefully support by the 
DFG priority program 'Global Methods in Complex Geometry'.}
\date{\today}

\begin{document}

\begin{abstract}
Tian and Yau constructed in \cite{ty} a complete Ricci-flat K\"ahler metric on the complement of an ample and smooth anticanonical divisor. We inquire
into the behaviour of this metric towards the boundary divisor and prove a slow decay rate of the difference to an appropriate explicitely given
referential metric.
\end{abstract}

\maketitle

\section{Introduction}

The problem addressed in this paper is connected to the problem of finding complete Ricci-flat K\"ahler metrics on open complex manifolds. 
As a model for this situation serves the complement of a divisor on a compact manifold $X$. If
$D\in|-K_X|$, then the section of $-K_X$ vanishing exactly on $D$ yields an isomorphism 
$\Omega^n_{X\setminus D}\cong{\mathcal O}_{X\setminus D}$. In analogy to the Calabi conjecture on compact 
manifolds this raises
the expectation that there exists a complete Ricci-flat K\"ahler metric on $X\setminus D$. Moreover, methods
to find such a metric may also work, if $D$ is not reduced, leading to the speculative existence of
Ricci-flat K\"ahler metrics on $X\setminus D$, whenever $-(D+K_X)$ is effective. 

Part of this program has already been established. Tian and Yau have proved in \cite{ty} and \cite{ty2}
the existence of a complete Ricci-flat K\"ahler metric in case $D\in|-K_X|$ is neat, almost ample and smooth.
Bando and Kobayashi have shown the claim in \cite{bk}, if $rD\in|-K_X|$ for $r>1$ and $D$ is ample, smooth and admits
on itself a K\"ahler-Einstein metric.

Of course, the techniques introduced in \cite{ty}, \cite{ty2} and \cite{bk} cannot be easily 
generalized to the reducible or non-reduced case.
The case of a reducible $D$ is still open. One
of the difficulties to treat the reducible case is the lack of a more detailed knowledge of the asymptotic behaviour towards $D$ of the constructed Ricci-flat metric 
in the smooth case. 
In this paper we will prove that the metric constructed by Tian and Yau obeys a polynomial
approach rate to an explicitely given initial metric: 
Let $D$ be given by the section $S\in H^0(-K_X)$. For an appropriately chosen metric $\|.\|$ on $-K_X$, the initial metric 
$\omega:=i\frac n{n+1}\dd(\log\|S\|^2)^{\frac{n+1}n}$ is complete and the complete Ricci-flat metric $\tilde\omega$ constructed in \cite{ty} satisfies

\begin{thm}$(1-C(-\log\|S\|^2)^{-\frac 1{6n}})\omega\le\tilde\omega\le (1+C(-\log\|S\|^2)^{-\frac 1{6n}})\omega$.\end{thm}

Like in \cite{ty} we consider the Monge-Amp\`ere equation
\begin{eqnarray*}
\frac{(\omega+i\partial\overline{\partial}u)^n}
{\omega^n}=e^{f}
\end{eqnarray*} 
for a fast decaying $f$.  

The first step to prove this is to use an appropriately disturbed equation that allows for a maximum principle.
This is given by
\begin{eqnarray*}
\frac{(\omega+i\partial\overline{\partial}u_{\varepsilon})^n}
{\omega^n}=e^{\varepsilon u_{\varepsilon}+f}
\end{eqnarray*} 
We prove strong decay for the unique bounded solution constructed in \cite{cy}
and give a description of the behaviour of the involved constants in dependence of $\varepsilon$. Here we make use of a weak maximum principle. 
By curvature estimates provided by \cite{ty} and again by a weak maximum principle we are able to prove the theorem.

The decay exponent obtained here is $\frac 1{6n}$ and hence too small in order to be useful for the case of a reducible divisor. Nevertheless
it is a first decay result. We do not believe that it is optimal, so we hope to be able to improve the decay in the future.   
The authors would like to thank Hans-Christoph Grunau for fruitful discussions and encouragement and Georg Schumacher for supporting the project.

\section{Construction of the initial metric}

In this section we use a slightly different approach for the construction of the initial metric used in \cite{ty}.
Let $X$ be smooth projective, $\dim_{\mathbb{C}}(X)=n$,
$D\subset X$ a smooth ample divisor and
let the adjoint bundle $K_X+[D]\cong\mathbb{C}$ be trivial. Let
$\|\quad\|$ be a smooth metric for the divisor bundle $[D]$ such that
$\omega(0)=i\partial\overline{\partial}\log\|\quad\|>0$ is positive
definite on all $X$. Let $D=V(S(z)=0)$ be a defining section.
Because $K_X+[D]\cong\mathbb{C}$ is trivial, there is up to scalars
a unique meromorphic $n$-form $\Omega=\frac{\Psi}{S}$
on $X$ with poles along $D$ of first order. Let $U\subset X$ be
a chart with $S(z)=z_1$ and $\Omega$ represented on $U$ by
\begin{eqnarray*}
\Omega=\frac{\Psi(z_1,\cdots,z_n)dz_1\wedge\cdots\wedge dz_n}{z_1}.
\end{eqnarray*} 
Then
\begin{eqnarray*}
\hat{\Omega}=\Psi(0,z_2,\cdots,z_n)dz_2\wedge\cdots\wedge dz_n
\end{eqnarray*}
produces a holomorphic $(n-1)$-form on $D$, which is everywhere
nonzero on $D$, because $K_D$ is trivial. After multiplying by a 
constant $a>0$ we can assume:
\begin{eqnarray*}
\int_D\left(\omega(0)|_D\right)^{n-1}=\int_D\hat{\Omega}\wedge
\overline{\hat{\Omega}}.
\end{eqnarray*}  
So $\Omega$ is now uniquely determined. $\omega(0)|_D$ produces
a K\"ahler class on $D$, so by the Theorem of Yau \cite{y78} there is a 
$C^{\infty}$-function $\varphi_0:D\longrightarrow\mathbb{R}$
with $\omega(0)|_D+i\partial\overline{\partial}\varphi_0>0$ and
\begin{eqnarray*}
\left(\omega(0)|_D+i\partial\overline{\partial}\varphi_0\right)^{n-1}=
\hat{\Omega}\wedge\overline{\hat{\Omega}}.
\end{eqnarray*}
Now according to a result of Schumacher \cite[Thm 4]{sch} the metric
$\|\quad\|_De^{-\varphi_0}$ for the restricted bundle $[D]|_D$
has an extension $\|\quad\|e^{-\varphi}$, whose curvature is positive on all of $X$ 
because $[D]$ is ample. So we replace the metric $\|\quad\|$
by the new bundle metric $\|\quad\|e^{-\varphi}$ and set $\omega(0)>0$ as the curvature
form of the new bundle metric. Then we have 
\begin{eqnarray*}
\omega(0)|_D^{n-1}=
\hat{\Omega}\wedge\overline{\hat{\Omega}}.
\end{eqnarray*}
Following \cite{ty} we define a complete K\"ahler metric by
$$\omega:=i\frac n{n+1}\dd(-\log\|S\|^2)^{\frac{n+1}{n}}.$$

First we observe
\begin{eqnarray*}
\omega&=&i\frac{n}{n+1}\partial\overline{\partial}
\left(-\log\|S\|^2\right)^{1+\frac{1}{n}}\\
&=&\omega(0)\left(-\log\|S\|^2\right)^{\frac{1}{n}}+
\frac{i}{n}\partial\log\|S\|^2\wedge
\overline{\partial}\log\|S\|^2
\left(-\log\|S\|^2\right)^{\frac{1}{n}-1}\\
&>&0.
\end{eqnarray*}
So $\omega$ defines a complete K\"ahler form.

We compute the volume form:
\begin{eqnarray*}
\omega^n=(\omega(0))^n\left(-\log\|S\|^2\right)+
(\omega(0))^{n-1}\wedge\partial\log\|S\|^2\wedge
\overline{\partial}\log\|S\|^2.
\end{eqnarray*}
Expanding in a coordinate chart $U$ with $S(z)=z_1$ we recognize
\begin{eqnarray*}
\omega^n=\frac{1}{|z_1|^2}
\left(\hat{\Omega}\wedge\overline{\hat{\Omega}}\wedge dz_1\wedge d\overline z_1+(z_1h_1+
\overline{z_1h_1})dz_1\wedge\dots\wedge dz_n\wedge d\overline z_1\wedge\dots\wedge d\overline z_n\right)
\end{eqnarray*}
with $\|h_1\|_{C^3_\omega(U)}\leq C$ bounded. Inasmuch we have 
\begin{eqnarray*}
\Omega\wedge\overline{\Omega}=\frac{1}{|z_1|^2}
\left(\hat{\Omega}\wedge\overline{\hat{\Omega}}\wedge dz_1\wedge d\overline z_1+(z_1H_1+
\overline{z_1H_1}))dz_1\wedge\dots\wedge dz_n\wedge d\overline z_1\wedge\dots\wedge d\overline z_n\right)
\end{eqnarray*}
with $\|H_1\|_{C^3_\omega(U)}\leq C$ bounded. Defining
\begin{eqnarray*}
F=\frac{\Omega\wedge\overline{\Omega}}{\omega^n}
\end{eqnarray*}
we see that $F$ has an expansion $F=1+z_1F_1+\overline{z_1F_1}$
with $\|F_1\|_{C^3_\omega(U)}\leq C$ bounded. This implies
that 
\begin{eqnarray*}
\int_X\left|\omega^n-\Omega\wedge\overline{\Omega}\right|<\infty.
\end{eqnarray*}
Retracing the construction carefully we see that by adding a constant to $\varphi$ we can achieve  
\begin{eqnarray}
\label{integral}\int_X\left(\omega^n-\Omega\wedge\overline{\Omega}\right)=0.
\end{eqnarray}

But now we have to ensure that
$\|S(z)\|\leq\frac{1}{e}$ on all of $X$. This we do by manipulating $\varphi$ again. This time we replace
$\varphi$ by $\varphi+C\|S\|^4$. This is an admissible extension of $\varphi_0$.
By an appropriate choice of $C\gg 0$ we can achieve $\|S(z)\|\le\frac 1e$ on all of $X$. Using the formula for any $\alpha=\beta+\gamma$, 
\begin{eqnarray*}\left(\frac{n}{n+1}\right)^n\left( (\partial\overline\partial\alpha^{\frac{n+1}n})^n-(\partial\overline\partial\beta^{\frac{n+1}n})^n\right)&=&
d\left(\sum_{k=1}^n{n\choose k}\beta\overline\partial\gamma\wedge(\dd\gamma)^{k-1}\wedge(\dd\beta)^{n-k}+\right.\\
& &\left. +{n\choose{k-1}}\gamma\overline\partial\beta\wedge(\dd\gamma)^{k-1}\wedge
(\dd\beta)^{n-k}\right.+\\
& & \left. +{n-1\choose k}\partial\beta\wedge\overline\partial\beta\wedge\overline\partial\gamma\wedge(\dd\gamma)^{k-1}\wedge(\dd\beta)^{n-k-1}+\right.\\
& &\left. +{{n-1}\choose{k-2}}\partial\gamma\wedge\overline\partial\gamma\wedge\overline\partial\beta\wedge(\dd\gamma)^{k-2}\wedge(\dd\beta)^{n-k}\right),\end{eqnarray*}
we obtain
for $\alpha:=-\log\|S\|^2_{\varphi}+C\|S\|^4$, 
$$(\partial\overline\partial\alpha^{\frac{n+1}n})^n=(\partial\overline\partial(-\log\|S\|^2_\varphi)^{\frac{n+1}n})^n+d(\|S\|^2\eta)$$
for an $(n-1,n)$-form $\eta$ with at most logarithmic singularities along $D$. By Stokes' Theorem, the integral condition (\ref{integral}) still holds. Finally, a similar calculation shows that $\omega^n$ and $\Omega\wedge\overline\Omega$
are still asymptotically equivalent.

Due to \cite{ty}, by changing $\varphi$ by higher powers of $\|S\|^2$, all the above mentioned properties can be preserved and the additional property
$$\log F\sim (-\log\|S\|^2)^{-N}$$
be obtained for arbitrarily big chosen $N$.

\section{Fast decay of the $\varepsilon$-solution}

Let
$f:=\log F$. By Cheng and Yau \cite{cy} the Monge-Amp\`ere equation
\begin{eqnarray*}
\frac{(\omega+i\partial\overline{\partial}u_{\varepsilon})^n}
{\omega^n}=e^{\varepsilon u_{\varepsilon}+f}
\end{eqnarray*} 
has a unique bounded $C^{\infty}$-solution
$u_{\varepsilon}:X\backslash D\longrightarrow\mathbb{R}$.

\begin{thm}\label{ueps}Let $u_\varepsilon$ be the bounded solution of 
$$\frac{\det(\omega+i\dd u)}{\det\omega}=e^{f+\varepsilon u}.$$
Then for all $0<\delta<N-1-\frac 1n$ there exists $C>0$ independent of $\varepsilon$ such that 
$$\max\{|u_\varepsilon|,\|du_\varepsilon\|_\omega, \|\dd u_\varepsilon\|_\omega\}\le C\varepsilon^{-(1+n\delta-\frac{\delta}{n+1})}(-\log\|S\|^2)^{-\delta}.$$
\end{thm}

\begin{proof} Consider for a $\beta>0$ the function
$v_{\beta\varepsilon}$ defined by
\begin{eqnarray*}
u_{\varepsilon}=v_{\beta\varepsilon}\left(\beta+
(-\log\|S\|^2)^{-\delta}\right).
\end{eqnarray*} 
\begin{prop}
The map
\begin{eqnarray*}
z\longmapsto\left(\beta+(-\log\|S\|^2)^{-\delta}\right)^{-1}
\end{eqnarray*}
is of bounded geometry.
\end{prop}
\begin{proof} Quasicoordinates are given for $C\gg 1$ and
$|w_1|\leq\pi$ by
\begin{eqnarray*}
w_1\longmapsto e^{-((w_1+C)^{\frac{2n}{n+1}})}=z_1.
\end{eqnarray*}
So if we denote by $h(w_1)$ the pull back of the above function
we obtain
\begin{eqnarray*}
h(w_1)\cong\left(\beta+
\left(2Re\left((w_1+C)^{\frac{2n}{n+1}}\right)
\right)^{-\delta}\right)^{-1}.
\end{eqnarray*}
For $C$ large enough we see that 
\begin{eqnarray*}
\frac{1}{2\beta}\leq|h(w_1)|\leq\frac{3}{2\beta}
\end{eqnarray*}
is bounded independent of $C$. Furthermore a calculation shows that
we have bounds for all derivatives depending on $\beta$ but not on
$C\gg 1$ and $w_1$.\end{proof}

We conclude that $v_{\beta\varepsilon}$ is as a product of two 
functions of bounded geometry is also within this class and the
weak maximum principle is applicable. We know that
$|f|\leq C(-\log\|S\|^2)^{-N}$ with $N\gg 2$. As a subsequent
discussion shows we can assume that the $C^0$-bounded function
$v_{\beta\varepsilon}$ attains its supremum as a maximum in the 
interior of $X\backslash D$ in a point $z_{\beta\varepsilon}$ and
furthermore that $v_{\beta\varepsilon}(z_{\beta\varepsilon})\geq 0$.
Then we have 
$(\partial v_{\beta\varepsilon})(z_{\beta\varepsilon})= 
(\overline{\partial} v_{\beta\varepsilon})(z_{\beta\varepsilon})=0$
and $(\partial\overline{\partial}
v_{\beta\varepsilon})(z_{\beta\varepsilon})\leq 0$ negative 
semidefinite. We apply
\begin{prop}
For symmetric $n\times n$-matrices $A,B$ with $A>0$, $A+B>0$
and $B\leq 0$ we have $\det(A+B)\leq\det A$.
\end{prop}

We evaluate the Monge-Amp\`ere equation 
\begin{eqnarray*}
\frac{(\omega+i\partial\overline{\partial}u_{\varepsilon})^n}
{\omega^n}=e^{\varepsilon u_{\varepsilon}+f}
\end{eqnarray*} 
in the point $z_{\beta\varepsilon}$ and yield in this point
\begin{eqnarray*}
\frac{(\omega+iv_{\beta\varepsilon}
\partial\overline{\partial}(-\log\|S\|^2)^{-\delta})^n}
{\omega^n}&\geq&\frac{(\omega+i\partial\overline{\partial}u_{\varepsilon})^n}
{\omega^n}\\
&=&e^{\varepsilon(\beta+(-\log\|S\|^2)^{-\delta})v_{\beta\varepsilon}+f}\\
&\geq&
1+\varepsilon(\beta+(-\log\|S\|^2)^{-\delta})v_{\beta\varepsilon}
-\tilde{f},
\end{eqnarray*}
where $|\tilde{f}|\leq C(-\log\|S\|^2)^{-N}$. As before we 
calculate
\begin{eqnarray*}
i\partial\overline{\partial}(-\log\|S\|^2)^{-\delta}&=&
-\delta\omega(0)(-\log\|S\|^2)^{-\delta-1}+\\
& &+i\delta(\delta+1)\partial\log\|S\|^2\wedge
\overline{\partial}\log\|S\|^2
(-\log\|S\|^2)^{-\delta-2}.
\end{eqnarray*}
We realize
\begin{eqnarray*}
\Delta_{\omega}\left((-\log\|S\|^2)^{-\delta}\right)=
(-\log\|S\|^2)^{-\delta-1-\frac{1}{n}}\delta(n\delta+1)+
O((-\log\|S\|^2)^{-\delta-2-\frac{1}{n}}
\end{eqnarray*}
and
\begin{eqnarray*}
\left|\frac{\omega^{n-k}\wedge\left(
\partial\overline{\partial}(-\log\|S\|^2)^{-\delta}\right)^k}
{\omega^n}\right|\leq
C(-\log\|S\|^2)^{-k(\delta+1)-\frac{k}{n}}.
\end{eqnarray*}
We have the trivial bound $|u_{\varepsilon}|\leq
\varepsilon^{-1}\max|f|=:C{\varepsilon}^{-1}$. Then
$|v_{\beta\varepsilon}|(-\log\|S\|^2)^{-\delta}\leq
\frac C{\varepsilon}$ is also bounded independent of $\beta>0$. We end up
with
\begin{eqnarray*}
\frac{(\omega+iv_{\beta\varepsilon}
\partial\overline{\partial}(-\log\|S\|^2)^{-\delta})^n}
{\omega^n}=1+v_{\beta\varepsilon}\delta(n\delta+1)
(-\log\|S\|^2)^{-\delta-1-\frac{1}{n}}\left(1+
V_{\beta\varepsilon}(z)\right)
\end{eqnarray*}
where $|V_{\beta\varepsilon}(z)|\leq\frac{C^{n-1}}{\varepsilon^{n-1}}
(-\log\|S\|^2)^{-1}$ and this equality holds for all $z\in
X\backslash D$. We combine the inequality coming from the 
Monge-Amp\`ere equation and the above equation to get (in the point
$z_{\beta\varepsilon}$):
\begin{eqnarray*}
\tilde{f}(z_{\beta\varepsilon})&\geq& 
v_{\beta\varepsilon}(z_{\beta\varepsilon})(\varepsilon
\left(\beta+(-\log\|S\|^2)^{-\delta}\right)-
\delta(n\delta+1)
(-\log\|S\|^2)^{-\delta-1-\frac{1}{n}}
(1+V_{\beta\varepsilon}(z_{\beta\varepsilon}))).
\end{eqnarray*}
There are two cases possible:\\
In the first case we have (in the point
$z_{\beta\varepsilon}$): 
\begin{eqnarray*}
\varepsilon
\left(\beta+(-\log\|S\|^2)^{-\delta}\right)\geq 
2\delta(n\delta+1)
(-\log\|S\|^2)^{-\delta-1-\frac{1}{n}}
(1+\frac{C^{n-1}}{\varepsilon^{n-1}}(-\log\|S\|^2)^{-1}).
\end{eqnarray*}
Then we conclude
\begin{eqnarray*}
v_{\beta\varepsilon}(z_{\beta\varepsilon})&\leq&
\tilde{f}(z_{\beta\varepsilon})\frac{1}{\delta(n\delta+1)}
(-\log\|S\|^2(z_{\beta\varepsilon}))^{\delta+1+\frac{1}{n}}
(1+\frac{C^{n-1}}{\varepsilon^{n-1}}(-\log\|S\|^2)^{-1})^{-1}\\
&\leq&\tilde{C}(1+\varepsilon^{n-1})
\end{eqnarray*}
because $|\tilde{f}(z_{\beta\varepsilon})|\leq
(-\log\|S\|^2(z_{\beta\varepsilon}))^{-N}$ for $0<\delta\leq 
N-1-\frac{1}{n}$. So in this case we have an upper bound 
for $v_{\beta\varepsilon}$ independent of $\beta>0$ satisfying even a better $\varepsilon$-behaviour than claimed.\\

In the complementary second case we get a lower bound for 
$(-\log\|S\|^2(z_{\beta\varepsilon}))^{-1}$:
\begin{eqnarray*}
\frac{1}{C^{n-1}}\varepsilon^{n-\frac{1}{n+1}}\leq
(-\log\|S\|^2(z_{\beta\varepsilon}))^{-1}.
\end{eqnarray*}
Now we have $u_{\varepsilon}(z_{\beta\varepsilon})\leq
\frac C{\varepsilon}$ and $v_{\beta\varepsilon}(z)\leq
v_{\beta\varepsilon}(z_{\beta\varepsilon})$.
This implies
\begin{eqnarray*}
u_{\varepsilon}(z)\left(\beta+(-\log\|S\|^2(z))^{-\delta}
\right)^{-1}\leq\frac C{\varepsilon}
\left(\beta+(-\log\|S\|^2
(z_{\beta\varepsilon}))^{-\delta}\right)^{-1}
\end{eqnarray*}
respectively
\begin{eqnarray*}
u_{\varepsilon}(z)\leq\frac C{\varepsilon}
\frac{\beta+(-\log\|S\|^2(z))^{-\delta}}
{\beta+(-\log\|S\|^2
(z_{\beta\varepsilon}))^{-\delta}}<
\frac C{\varepsilon}
\frac{\beta+(-\log\|S\|^2(z))^{-\delta}}
{\beta+\left(\frac{1}{C^{n-1}}\varepsilon^{n-\frac{1}{n+1}}
\right)^{\delta}}.
\end{eqnarray*}
Now we let $\beta\rightarrow 0$ tend to 0 and obtain
\begin{eqnarray*}
u_{\varepsilon}(z)\leq C^{1+\delta(n-1)}
\varepsilon^{-(1+n\delta-\frac\delta{n+1})}(-\log\|S\|^2(z))^{-\delta}.
\end{eqnarray*}
The same line of thought applies to the minimum of 
$v_{\beta\varepsilon}(z)$ so we get also lower bounds.
The estimates of the derivatives are provided by Schauder theory
in the local quasicoordinates as explained in \cite[Section 6]{sch}.\\

Now we want to discuss the case that $v_{\beta\varepsilon}(z)$ has
no maximum in $X\backslash D$ but takes its supremum on $D$.
So we assume $0<\sup\{v_{\beta\varepsilon}(z)\}=L<\infty$ and 
$v_{\beta\varepsilon}(z)<L$ for all $z\in X\backslash D$.
According to Cheng and Yau \cite{cy} we have a sequence of points 
$z(m)\in X\backslash D$, $m\in\mathbb{N}$ with the following
properties: For $m\seq\infty$
\begin{enumerate}
\item $v_{\beta\varepsilon}(z(m))\rightarrow L$,
\item $(\partial v_{\beta\varepsilon}(z(m))\rightarrow 0\mbox{ and }
(\overline{\partial} v_{\beta\varepsilon}(z(m))\rightarrow 0$,
\item $\overline{\lim}\partial\overline{\partial}v_{\beta\varepsilon}(z(m))\leq
0\mbox{ is negative semidefinite}$.
\end{enumerate}
We put $\check{\omega}=\omega+\partial\overline{\partial}
u_{\varepsilon}$ and $\hat{\omega}=\check{\omega}-
v_{\beta\varepsilon}\partial\overline{\partial}
(-\log\|S\|^2(z))^{-\delta}$. By Cheng and Yau \cite{cy} we have
$C^2$-estimates $\frac{1}{\hat{C}_{\varepsilon}}\omega\leq
\check{\omega}\leq\hat{C}_{\varepsilon}\omega$. Furthermore from
$|v_{\beta\varepsilon}|(-\log\|S\|^2(z))^{-\delta}\leq 
C_{\varepsilon}$ and from the asymptotics of 
$\partial\overline{\partial}(-\log\|S\|^2(z))^{-\delta}$ we
conclude an estimate $\frac{1}{\tilde{C}_{\varepsilon}}\omega\leq
\hat{\omega}\leq\tilde{C}_{\varepsilon}\omega$. We develop as before
\begin{eqnarray*}
\check{\omega}^n&=&\left(\hat{\omega}+iv_{\beta\varepsilon}
\partial\overline{\partial}(-\log\|S\|^2(z))^{-\delta}\right)^n\\
&=&\hat{\omega}^n+\sum_{j=1}^n
\left(
\begin{array}{c}
n\\
j
\end{array}
\right)\hat{\omega}^{n-j}\wedge v_{\beta\varepsilon}^j
\left(i\partial\overline{\partial}(-\log\|S\|^2(z))^{-\delta}\right)^j.
\end{eqnarray*}
Looking at the asymptotics we see
\begin{eqnarray*}
\left|\frac{1}{\omega^n}\sum_{j=1}^n
\left(
\begin{array}{c}
n\\
j
\end{array}
\right)\hat{\omega}^{n-j}\wedge v_{\beta\varepsilon}^j
\left(\partial\overline{\partial}(-\log\|S\|^2(z))^{-\delta}\right)^j
\right|\leq\tilde{C}_{\varepsilon}
(-\log\|S\|^2(z))^{-\delta-1-\frac{1}{n}}.
\end{eqnarray*}
For $m\geq m_0$ large we have $z(m)\in U_{\varepsilon}(D)$ and
$v_{\beta\varepsilon}(z(m))>0$. Then the Monge-Amp\`ere equation implies
\begin{eqnarray*}
\frac{\hat{\omega}^n}{\omega^n}(z(m))&\geq& 1-|f(z(m))|+
v_{\beta\varepsilon}(z(m))\cdot\\
& &\left(\varepsilon(\beta+
(-\log\|S\|^2(z))^{-\delta})-\tilde{C}_{\varepsilon}
(-\log\|S\|^2(z))^{-\delta-1-\frac{1}{n}}\right)
\end{eqnarray*}
or, respectively,
\begin{eqnarray*}
\left(\varepsilon(\beta+
(-\log\|S\|^2(z))^{-\delta})-\tilde{C}_{\varepsilon}
(-\log\|S\|^2(z))^{-\delta-1-\frac{1}{n}}\right)^{-1}
\left(\frac{\hat{\omega}^n}{\omega^n}(z(m))-1\right)
+& &\\
|f(z(m))|\left(\varepsilon(\beta+
(-\log\|S\|^2(z))^{-\delta})-\tilde{C}_{\varepsilon}
(-\log\|S\|^2(z))^{-\delta-1-\frac{1}{n}}\right)^{-1}&\geq&
v_{\beta\varepsilon}(z(m))
\end{eqnarray*}
Now we hold $\varepsilon$,$\beta>0$ fix and take the limit 
$m\rightarrow\infty$:
\begin{eqnarray*}
v_{\beta\varepsilon}(z(m))&\rightarrow& L>0\\
\varepsilon(\beta+
(-\log\|S\|^2(z))^{-\delta})-\tilde{C}_{\varepsilon}
(-\log\|S\|^2(z))^{-\delta-1-\frac{1}{n}}&\rightarrow&
\varepsilon\beta\\
|f(z(m))|&\rightarrow& 0\\
\overline{\lim}\left(\frac{\hat{\omega}^n}{\omega^n}(z(m))-1\right)&\leq& 0.
\end{eqnarray*}
So we conclude:
\begin{eqnarray*}
0\geq
\overline{\lim}\left(\frac{\hat{\omega}^n}{\omega^n}(z(m))-1\right)\geq
\varepsilon\beta L>0,
\end{eqnarray*}
which is a contradiction, so there must be a maximum inside of
$X\backslash D$.
\end{proof}

\section{Slow decay of the solution}

Now we study the original Monge-Amp\`ere equation for our problem
$$\frac{(\omega+i\dd u)^n}{\omega^n}=e^f.$$
By \cite{ty} there exists a unique bounded $C^2$-solution of the problem such that $\tilde\omega:=\omega+i\dd u>0$. 
We want to study its decay to $D$. 

Unfortunately, the estimate of Proposition \ref{ueps} is not uniform in $\varepsilon$, so we have to add arguments in order to achieve our result. 
Our result here is much weaker than the properties of $u_\varepsilon$ provided by Theorem \ref{ueps}. We begin
with the preparation of elementary tools. 

\begin{lemma}\label{sqrteps}There are constants $C_n>0$ such that for all positive real numbers $a_1,\dots,a_n$ and $0<\varepsilon<1$ satisfying 
$$\sum_{i=1}^na_i\le n(1+\varepsilon), \prod_{i=1}^na_i=1,$$
holds: $(1+C_n\sqrt\varepsilon)^{-1}\le a_i\le 1+C_n\sqrt\varepsilon$ for every $i$.
\end{lemma}

\begin{proof}
This is easily seen by the inequality
$$\sum_{i=1}^n(\sqrt{a_i}-1)^2+2(\sum_{i=1}^n\sqrt{a_i}-n(\prod_{j=1}^na_j)^{\frac 1{2n}})\le n\varepsilon.$$
\end{proof}

As an immediate consequence we obtain that decay of $\Delta_\omega u$ is enough in order to have decay of $\|\dd u\|_\omega$:

\begin{lemma}\label{laplace}Let $C>0$ and $0<\beta<N$. If $\Delta_\omega u\le C(-\log\|S\|^2)^{-\beta}$, then
$$(1-C(-\log\|S\|^2)^{-\frac\beta 2})\omega\le\omega+\dd u\le (1+C(-\log\|S\|^2)^{-\frac\beta 2})\omega.$$
\end{lemma}

\begin{proof}
Let us choose a point $x\in X$ and coordinates such that in $\omega(x)=\delta_{ij}, \dd u(x)=u_{ij}\delta_{ij}$. Let us denote
$a_i:=1+u_{,i\bar i}$. By $\omega+\dd u>0$ we know $a_i>0$. Furthermore,
$$\prod a_i=e^f=1+O((-\log\|S\|^2)^{-N}), \sum a_i=n+\Delta_\omega u\le n+C(-\log\|S\|^2)^{-\beta}.$$
So we can change $a_i$ in such a way that $b_i:=a_i(1+O((-\log\|S\|^2)^{-N}))$ satisfies $\prod b_i=1$. Of course, still
$\sum b_i\le n+C(-\log\|S\|^2)^{-\beta}$. Now we apply Lemma \ref{sqrteps} and obtain
$$u_{,i\bar i}=O((-\log\|S\|^2)^{-\frac\beta 2}).$$
This implies the claim.
\end{proof}

\begin{lemma}\label{aij}
If $A\in M(n,{\bb C})$ is a hermitian matrix and $\|A\|:=\sup_{x\in{\bb C}^n}\frac{|(x,Ax)|}{\|x\|^2}$, then the estimate holds $|a_{ij}|\le n\|A\|$ for
any $i,j$.
\end{lemma}

\begin{thm}$(1-C(-\log\|S\|^2)^{-\frac 1{6n}})\omega\le\omega+\dd u\le (1+C(-\log\|S\|^2)^{-\frac 1{6n}})\omega$.\end{thm}

\begin{proof}
Recall that $u$ is the bounded solution of
$$\frac{\det(\omega+i\dd u)}{\det\omega}=e^f$$
and the curvature tensor is given by
$$R_{i\bar j k\bar l}=\partial_i\overline\partial_j\omega_{kl}-\omega^{pq}\partial_i\omega_{kq}\overline\partial_j\omega_{pl}.$$
By a computation of \cite{ty},
$$\|R^{k\bar l}_{i\bar j}\|_\omega=\|R\4 ijkl\|_\omega\le C(-\log\|S\|^2)^{-\frac 1n}.$$
Here $\|R\4 ijkl\|_\omega$ is the norm of the bisectional curvature $\sup\frac{|R\4 ijkl\zeta^i\overline\zeta^j\xi^k\overline\xi^l|}{g(\zeta,\zeta)g(\xi,\xi)}$. 
Now we choose coordinates centered around $x$ such that $\omega_{ij}(x)=\delta_{ij}, \omega_{ij,k}(x)=0$. For $\xi=e_k$ the estimate of the bisectional curvature
yields
$$\|(R\4 ijkk)_{ij}\|\le C(-\log\|S\|^2)^{-\frac 1n}$$
in the sense of Lemma \ref{aij}. The matrix $(\sqrt{-1}R\4 ijkk)_{ij}$ is hermitian, indeed, hence
\begin{equation}\label{rijkk}|R\4 ijkk|\le C(-\log\|S\|^2)^{-\frac 1n}\end{equation}
for any $i,j,k$. 

Further, we obtain in $x$
\begin{eqnarray*}\Delta_{\tilde\omega}\Delta_\omega u&=&\tilde\omega^{ij}\partial_i\overline\partial_j(\omega^{kl}\partial_k\overline\partial_l u)\\
&=&\tilde\omega^{ij}\partial_i\overline\partial_j(\omega^{kl}(\tilde\omega_{kl}-\omega_{kl}))\\
&=&\tilde\omega^{ij}(\tilde\omega_{kl}-\omega_{kl})\partial_i\overline\partial_j\omega^{kl}+\tilde\omega^{ij}\omega^{kl}\partial_i\overline\partial_j\partial_k\overline
\partial_l u\\
&=&\tilde\omega^{ij}R^{k\bar l}_{i\bar j}\tilde\omega_{kl}-\tilde\omega^{ij}R_{ij}+\tilde\omega^{ij}\omega^{kl}\partial_i\overline\partial_j\partial_k\overline\partial_l u\\
\end{eqnarray*}
by the choice of the coordinates. Now we have to replace the fourth derivatives of $u$.
Since $\tilde\omega$ is Ricci-flat, we may use the equations
$$0=\tilde R\ii kl=\partial_k\overline\partial_l\log\det(\omega+i\dd u)=\partial_k\tilde\omega\II ij\overline\partial_l\partial_i\overline\partial_j u+
\tilde\omega\II ij(\partial_k\overline\partial_l\omega\ii ij+\partial_i\overline\partial_j\partial_k\overline\partial_l u).$$
Now we substitute $\partial_k\tilde\omega\II ij=-\tilde\omega\II im\tilde\omega\II nj\partial_k\overline\partial_m\partial_n u$ and obtain for any $k,l$
$$\tilde\omega\II ij\partial_i\overline\partial_j\partial_k\overline\partial_l u=
\tilde\omega\II im\tilde\omega\II nj\overline\partial_l\partial_i\overline\partial_j u\cdot\partial_k\overline\partial_m\partial_n u-\tilde\omega\II ijR\4 klij.$$
So,
$$\Delta_{\tilde\omega}\Delta_\omega u=\tilde\omega^{ij}R^{k\bar l}_{i\bar j}\tilde\omega_{kl}-\omega\II kl\tilde\omega\II ijR\4 klij
-\tilde\omega^{ij}R_{ij}+\omega\II kl\tilde\omega\II im\tilde\omega\II nj\overline\partial_l\partial_i\overline\partial_j u\cdot\partial_k\overline\partial_m\partial_n u.$$

From now on we specialize our coordinates further in such a way that $\partial_i\overline\partial_j u=\lambda_i\delta_{ij}$ in $x$; this is possible by a locally constant
unitary coordinate transformation. The $C^2$-estimate of \cite{ty} yield, that $0<c\le 1+\lambda_i\le C$ with constants $c,C$ independent of $x$. 
The first two terms are now easily dealt with by (\ref{rijkk}). Observe that $f$ is explicitely computable in our case. In particular
$$|R_{ij}|=|\partial_i\overline\partial_j f|\le C(-\log\|S\|^2)^{-N}.$$
So the only mysterious term is
\begin{eqnarray*}a&:=&\omega\II kl\tilde\omega\II im\tilde\omega\II nj\overline\partial_l\partial_i\overline\partial_j u\cdot\partial_k\overline\partial_m\partial_n u\\
&=&\sum_{i,j,k}(1+\lambda_i)^{-1}(1+\lambda_j)^{-1}|\overline\partial_k\partial_i\overline\partial_ju|^2\end{eqnarray*}

Fortunately, $a\ge 0$ and we will use this to our advantage in the computation
\begin{eqnarray*}\Delta_{\tilde\omega}((-\log\|S\|^2)^\delta\Delta_\omega u)&=&(-\log\|S\|^2)^\delta\Delta_{\tilde\omega}\Delta_\omega u
+\Delta_{\tilde\omega}(-\log\|S\|^2)^\delta\Delta_\omega u+\\
& &+2(\nabla_{\tilde\omega}(-\log\|S\|^2)^\delta,\nabla_{\tilde\omega}\Delta_\omega u)_{\tilde\omega}.\end{eqnarray*}
Like above, we have the estimates
$$|\Delta_{\tilde\omega}(-\log\|S\|^2)^\gamma|\le(-\log\|S\|^2)^{\gamma-1-\frac 1n}, 
\|\nabla_{\tilde\omega}(-\log\|S\|^2)^\gamma\|^2\le(-\log\|S\|^2)^{2\gamma-1-\frac 1n}.$$
In order to make use of Cauchy-Schwarz we compute
\begin{eqnarray*}\|\nabla_{\tilde\omega}\Delta_\omega u\|^2_{\tilde\omega}&=&\|d\Delta_\omega u\|^2_{\tilde\omega}\\
&=&\tilde\omega\II ij\partial_i\Delta_\omega u\overline\partial_j\Delta_\omega u\\
&=&\tilde\omega\II ij\omega\II kl\omega\II mn\partial_i\partial_k\overline\partial_lu\cdot\overline\partial_j\partial_m\overline\partial_nu\\
&=&\sum_{i,k,m}(1+\lambda_i)^{-1}\partial_i\partial_k\overline\partial_ku\cdot\partial_i\partial_m\overline\partial_m u\\
&\le&\frac 12\sum_{i,k,m}(1+\lambda_i)^{-1}(|\partial_i\partial_k\overline\partial_ku|^2+|\partial_i\partial_m\overline\partial_mu|^2)\\
&=&\sum_{i,k}(1+\lambda_i)^{-1}|\partial_i\partial_k\overline\partial_ku|^2\end{eqnarray*}
The coefficient of $|\partial_i\partial_k\overline\partial_ku|^2$ in $a$ is $(1+\lambda_k)^{-1}((1+\lambda_i)^{-1}+(1+\lambda_k)^{-1})$. Hence, for
some constant $D$ depending on $c$ and $C$ we obtain
$$\|\nabla_{\tilde\omega}\Delta_\omega u\|_{\tilde\omega}\le D\sqrt a.$$
Now this
yields together with Cauchy-Schwarz
\begin{eqnarray*}\Delta_{\tilde\omega}((-\log\|S\|^2)^\delta\Delta_\omega u)&\ge&(-\log\|S\|^2)^\delta a-(-\log\|S\|^2)^{\delta-\frac 1n}-(-\log\|S\|^2)^{-N}-\\
& &-(-\log\|S\|^2)^{\delta-1-\frac 1n}-(-\log\|S\|^2)^{\delta-\frac 12-\frac 1{2n}}\sqrt a.\end{eqnarray*}
The sum of $a$-terms has a minimum as a function of $a$. This is attained at $a=O((-\log\|S\|^2)^{-1-\frac 1n})$ and hence
$$\Delta_{\tilde\omega}((-\log\|S\|^2)^\delta\Delta_\omega u)\ge-(-\log\|S\|^2)^{\delta-\frac 1n}.$$
Since $\Delta_{\tilde\omega}u=\sum\frac{\lambda_i}{1+\lambda_i}=n-\sum\frac 1{1+\lambda_i}$,
$$\Delta_{\tilde\omega}((-\log\|S\|^2)^{\delta}\Delta_\omega u-u)\ge-(-\log\|S\|^2)^{\delta-\frac 1n}-n+\sum\frac 1{1+\lambda_i}.$$
Now we would like to choose $x$ as a maximum of $(-\log\|S\|^2)^{\delta}\Delta_\omega u-u$ but, of course, this would assume what we would like to prove.
So we do the same considerations for $u_\varepsilon$ instead of $u$. We denote $\omega_\varepsilon=\omega+i\dd u_\varepsilon$. 
All computations are the same except that $R^{(\varepsilon)}\ii ij\not=0$ now, but 
$$\Ric(\omega_\varepsilon)=\Ric\omega+i\dd(f+\varepsilon u_\varepsilon)=i\varepsilon\dd u_\varepsilon.$$
So, by \cite{ty} we get $\|\Ric(\omega_\varepsilon)\|_\omega\le C\varepsilon$. On the other hand,
by Proposition \ref{ueps} we know $\|\Ric\omega_\varepsilon\|_\omega\le C\varepsilon^{-(n\alpha-\frac{\alpha}{n+1})}(-\log\|S\|^2)^{-\alpha}$. 
We combine both inequalities to obtain
$$\|\Ric(\omega_\varepsilon)\|_\omega\le C(-\log\|S\|^2)^{-\frac{\alpha}{1+\alpha(n-\frac 1{n+1})}}\le C(-\log\|S\|^2)^{-\frac 1n}$$
for $\alpha\ge n+1$ and with a constant $C$ independent of $\varepsilon$.
Hence we obtain, in particular,
$$\Delta_{\tilde\omega}((-\log\|S\|^2)^{\delta}\Delta_\omega u_\varepsilon-u_\varepsilon)\ge-C(-\log\|S\|^2)^{\delta-\frac 1n}-n+\sum\frac 1{1+\lambda^{(\varepsilon)}_i}$$
with $C$ independent of $\varepsilon$! 
By the decay property of $u_\varepsilon$ and $\Delta_\omega u_\varepsilon$ we can find a point $x_\varepsilon\in X$ where
$(-\log\|S\|^2)^{\delta}\Delta_\omega u_\varepsilon-u_\varepsilon$ attains its maximum.
In $x_\varepsilon$ we obtain
$$C(-\log\|S\|^2)^{\delta-\frac 1n}+n\ge \sum\frac 1{1+\lambda^{(\varepsilon)}_i}.$$
We know
$$\prod\frac 1{1+\lambda^{(\varepsilon)}_i}=e^{-(f+\varepsilon u_\varepsilon)}$$
and again we combine $|\varepsilon u_\varepsilon|\le C\varepsilon$ and 
$|\varepsilon u_\varepsilon|\le C\varepsilon^{-(n\alpha-\frac{\alpha}{n+1})}(-\log\|S\|^2)^{-\alpha}$ by an adequate multiplication to obtain
$$|\varepsilon u_\varepsilon|\le  C(-\log\|S\|^2)^{-\frac{\alpha}{1+\alpha(n-\frac 1{n+1})}}\le C(-\log\|S\|^2)^{-\frac 1n}$$
with $C$ independent of $\varepsilon$. 
Now by Lemma \ref{sqrteps} we obtain for 
$q:=Ce^{\frac 1n(f+\varepsilon u_\varepsilon)}(-\log\|S\|^2)^{\delta-\frac 1n}+e^{\frac 1n(f+\varepsilon u_\varepsilon)}-1$
\begin{equation*}\left(1+C'\sqrt{q}\right)^{-1}\le\frac 1{1+\lambda_i^{(\varepsilon)}}\le1+C'\sqrt{q}.\end{equation*}
So
$$\lambda_i^{(\varepsilon)}\le C'\sqrt{q}\le\tilde C(-\log\|S\|^2)^{\frac 12(\delta-\frac 1n)},$$
still only in the point $x_\varepsilon$.
In particular, if we assume $\delta=-\frac 12(\delta-\frac 1n)$, i.e. $\delta=\frac 1{3n}$ then for any $x\in X$
$$(-\log\|S\|^2)^{\delta}\Delta_\omega u_\varepsilon-u_\varepsilon\le C,$$
since $u_\varepsilon$ is uniformly bounded. We obtain also in the limit $\varepsilon\seq 0$
$$\Delta_\omega u\le C(-\log\|S\|^2)^{-\delta}.$$
Now we apply Lemma \ref{laplace} and obtain the desired result with decay rate $\frac 12\delta=\frac 1{6n}$. 
%
%
\end{proof}

\end{document}